\documentclass[12pt]{article}

\usepackage{amsmath}
\usepackage{amsthm}
\usepackage{amssymb}

\theoremstyle{plain}
\newtheorem{thm}{Theorem}[section]
\newtheorem{prop}[thm]{Proposition}
\newtheorem{cor}[thm]{Corollary}
\newtheorem{lem}[thm]{Lemma}

\theoremstyle{definition}
\newtheorem{defn}[thm]{Definition}
\newtheorem{ex}[thm]{Example}

\theoremstyle{remark}

\def\co{\colon\thinspace}

\title{Translation lengths in $Out(F_n)$}
\date{\today}
\author{Emina Alibegovi\'{c}}

\begin{document}
\pagebreak
\renewcommand{\thefootnote}{\null}
\maketitle
\setcounter{footnote}{0}
\renewcommand{\thefootnote}{\arabic{footnote}}

\begin{abstract}
We prove that all elements of infinite order in $Out(F_n)$ have positive 
translation lengths; moreover, they are bounded away from zero. Consequences 
include a new proof that solvable subgroups of $Out(F_n)$ are finitely 
generated and virtually abelian and the new result that such subgroups are 
quasi-convex. 
\end{abstract}

\section{Introduction} \label{intro}

In this paper we will study the translation lengths of outer automorphisms 
of a free group. Following \cite{gs3} we define the translation length 
$\tau_{X,G}(g)$ of $g \in \Gamma$ to be  
\[ \lim_{n \rightarrow \infty} \frac{\|g^n \|}{n} \]
where $\Gamma$ is a group with finite generating set $X$, and  
$\|g \|$ denotes the length of $g$ in the word metric on $\Gamma$ associated 
to $X$. 

Farb, Lubotzky and Minsky proved that Dehn twists (more generally, all 
elements of infinite order) in  $Mod(\Sigma_g)$ have positive translation 
length (\cite{flm}). We prove

\begin{thm}\label{thm1}
Every infinite order element $\mathcal{O} \in Out(F_n)$  has positive 
translation length. Furthermore, there exists a positive constant 
$\varepsilon_n$ such that $\tau(\mathcal{O}) \geq \varepsilon_n$. 

\end{thm}

Once more we can see the strong analogy between mapping class 
group of a surface, $Mod(\Sigma_g)$, and outer automorphism group of a free 
group, $Out(F_n)$.

To prove their theorem, Farb, Lubotzky and Minsky found a way to measure how 
much  a Dehn twist is `twisted' by looking at simple closed curves and 
their intersection number. Such an approach cannot work in the case 
of $Out(F_n)$ as we do not have an analogue of the intersection number. 

As a consequence of our main result we have 
\begin{cor}\label{cor1}
Every solvable subgroup of $Out(F_n)$ is finitely generated and 
virtually abelian.
\end{cor}
\noindent
Corollary~\ref{cor1} was proved in \cite{sol}, but Theorem~\ref{thm1} offers 
an alternative proof. 

\begin{cor}
Every abelian subgroup $A$ of $Out(F_n)$ is quasi-convex.
\end{cor}
\noindent
The proofs use techniques of \cite{gs3} and follow the same lines as the 
corresponding proofs in \cite{artin}.

I would like to thank Peter Brinkmann for suggesting  a careful examination 
of the exponents (see Definition~\ref{def1}). I also express gratitude to 
Mladen Bestvina for his support and help.  

\section{Translation lengths}

From the definition of translation length we can see that it depends on 
the choice of generating set for a group $\Gamma$. We will omit the 
reference to the generating set, since it will be clear which one we 
are using.

We list some properties of translation lengths which can be found in 
\cite{gs3}. 

\begin{prop}\label{prop1}
Let $X$ be a generating set for a group $\Gamma$. 
\begin{enumerate}
\item $0 \leq \tau(g) \leq \|g \|$
\item For all $x, g \in G$, $\tau(xgx^{-1})= \tau(g)$.
\item $\tau (g^n)= n\cdot \tau(g) \ \forall n \in \mathbb{N}$ . 
\end{enumerate}
\end{prop}

\noindent
Let $X=\{x_1, x_2, \ldots, x_n\}$ be a set of generators of a free group $F_n$.
Let $Y$ be the set of generators for $Aut(F_n)$ consisting of:
\begin{enumerate}
\item permutations ($x_i \mapsto x_j$, $x_j \mapsto x_i$, $x_k \mapsto x_k$ 
for all $k \neq i, j$), 
\item inversions ($x_i \mapsto x_i^{-1}$, $x_j \mapsto x_j$ for all $j \neq i$), 
\item Dehn twists ($x_i \mapsto x_i x_j$, $x_k \mapsto x_k$ for all $k \neq i$). 
\end{enumerate}
Let $\tilde{Y}$ denote the generating set for $Out(F_n)$ consisting of 
equivalence classes of elements of $Y$. 

Our goal is to prove that every element of infinite order in $Out(F_n)$ has 
positive translation length.  Since $Aut(F_n)$ embeds into $Out(F_{n+1})$, 
it will follow that every infinite order element of $Aut(F_n)$ has positive 
translation length. 

We will need the following definition for our proof:

\begin{defn}\label{def1}
Define a map $\alpha: F_n \rightarrow \mathbb{N}$ by
\begin{center}
\( \alpha (w)= \max \{ |p|: \tilde{w}^p\) is a subword of \( w \}\, , \)
\end{center}
where elments of $F_n$ are regarded as reduced words in the generators and 
their inverses. We also define 
\begin{center}
\( \tilde{\alpha}([w])= \max \{\alpha(u): u \) is a cyclically reduced 
conjugate of \( w \} \)
\end{center}
for the conjugacy class, $[w]$, of $w$. 
\end{defn}

\begin{lem}\label{lem1}
There exists a constant $C > 0$ such that for any $\tilde{g} \in\tilde{Y}$ 
and any cyclically reduced word $w \in F_n$ we have
\[ \tilde{\alpha}(\tilde{g}([w])) \leq \tilde{\alpha}([w]) + C\, . \]
\end{lem}

\begin{proof}
Note that inversions and permutations do not affect $\alpha(w)$, so we need 
only consider the case where $g$ is a Dehn twist.

Let $w \in F_n$ be a cyclically reduced element of length $n$. Let 
$w=A\, \tilde{w}^p\, B$, where $\alpha(w)=|p|$. Consider  
\[g(w)=[[ g(A)]]\, [[ g(\tilde{w})^p]]\, [[g(B)]]\, , \]
where $[[g(w)^p]]$ denotes the reduced word obtained from $g(w)^p$. 
By the \emph{Bounded Cancellation Lemma} (\cite{cooper}) there is a constant 
$C(g)$ such that at most $C(g)$ cancellations occur after concatenation of 
the words $[[ g(A)]]$ and $[[ g(\tilde{w})^p]]$. 
Hence $p$ can decrease by at most $2\, C(g)$ (cancellations may 
occur at the beginning and at the end of $[[ g(\tilde{w})^p]]$). Let 
$C_g= 2\, \max\{C(g), C(g^{-1})\}$. We now have 
\begin{align*}
&\alpha([[g(w)]]) \geq \alpha(w) - C_g \\ 
&\alpha(w) =\alpha(g^{-1}(g(w))) \geq \alpha([[g(w)]]) - C_g \\ 
&\alpha([[g(w)]]) \leq \alpha(w)+C_g. 
\end{align*}
If we take $\tilde{C}=\max \{C_g: g\  \text{a Dehn twist in}\  Y\}$, our 
claim is proved for elements of $Y$. 
Using a similar argument, we see that there is a constant $C$ such that
\[\tilde{\alpha}(\tilde{g}([w])) \leq \tilde{\alpha}([w])+ C\,. \] 
\end{proof}

\begin{ex}
We illustrate the idea of the proof of Theorem~\ref{thm1} with an example of 
a Dehn twist. Let $g$ be a Dehn twist which sends $x_2$ to $x_2x_1$ and 
fixes all other generators of $F_n$. 
$$\alpha(g^n(x_2))=\alpha(x_2x_1^n)=n\, .$$
If $g^n=g_1 \cdots g_m$, then $\|g^n\|=m$. By Lemma~\ref{lem1}, we have that 
\begin{align*}
& n=\alpha(g^n(x_2)) \leq \alpha(x_2)+m\, C=m\,C+1 \, ,\\
& \tau(g)=\lim_{n \rightarrow \infty} \frac{\|g^n\|}{n}\geq\lim_{n \rightarrow 
\infty}\frac{n-1}{n\, C}=\frac{1}{C}>0\, .
\end{align*}
So $g$ has positive translation length. 
\end{ex}

We give a short list of definitions which will be used throughout the rest
of the paper, but we suggest that the reader look at \cite{tits2}.   

Every element $\mathcal{O} \in Out(F_n)$ can be represented by a homotopy
equivalence $f \co G \rightarrow G$ of a graph $G$ whose fundamental group is 
identified with $F_n$. 
A map $\sigma \co J \rightarrow G$ ($J$ is an interval) is 
called a \emph{path} if it is either locally injective or a constant map 
(we also require that the endpoints of $\sigma$ are at vertices). 
Every map $\sigma \co J \rightarrow G$ is homotopic (relative endpoints) to a 
path $[[\sigma]]$. 

If $\sigma= \sigma_1 \ldots \sigma_l$ is a decomposition of a path or a 
circuit $\sigma$ into nontrivial subpaths we say that it is a 
\emph{k-splitting} if 
\[f^k(\sigma)=[[f^k(\sigma_1)]] \ldots [[f^k(\sigma_l)]]\]
is a decomposition into subpaths and is a \emph{splitting} if it is a 
$k$-splitting for all $k>0$. 

We say that a nontrivial path $\sigma \in G$ is a \emph{Nielsen path} for 
$f:G \rightarrow G$ if $[[f(\sigma)]]=\sigma$. The Nielsen path $\sigma$ is 
\emph{indivisible} if it cannot be written as a concatenation of nontrivial 
Nielsen paths. 

Let $\empty = G_0 \subsetneqq G_1 \subsetneqq \cdots \subsetneqq G_K = G$ be a 
filtration of $G$ by $f$-invariant subgraphs, and let $H_i=
\overline{G_i \backslash G_{i-1}}$.  Suppose $H_i$ is a single edge $E_i$ 
and $f(E_i)=E_i \upsilon^l$ for some closed indivisible Nielsen path 
$\upsilon \subset G_{i-1}$ and some $l>0$. The \emph{exceptional paths} 
are paths of the form $E_i \upsilon^k \overline{E_j}$ or 
$E_i \overline{\upsilon}^k \overline{E_j}$, where $k \geq0$, $j \leq i$ and 
$f(E_j)=E_j \upsilon^m$, for $m>0$.   
 
We remind the reader that every element of $Out(F_n)$ of infinite order has 
either exponential or polynomial growth (\cite{hb1}). A polynomially growing 
outer automorphism $\mathcal{O}\in Out(F_n)$ is unipotent if its action in 
$H_1(F_n;\mathbb{Z})$ is unipotent ($UPG$ automorphism). 

The following Theorem can be found in \cite{tits1}(page 564). 
\begin{thm}\label{thm2}
Suppose that $\mathcal{O} \in Out(F_n)$ is a $UPG$ automorphism. Then there is 
a topological representative $f \co G\rightarrow G$ of $\mathcal{O}$ with
the following properties: 
\begin{enumerate}
\item Each $G_i$ is the union of $G_{i-1}$ and a single edge $E_i$ satisfying
$f(E_i)=E_i\cdot u_i$ for some closed path $u_i$ that crosses only edges in 
$G_{i-1}$ ($\ \cdot$ indicates that the decomposition in question is a 
splitting). 
\item If $\sigma$ is any path with endpoints at vertices, then there exists 
$M=M(\sigma)$ so that for each $m\geq M$, $[[f^m(\sigma)]]$ splits into 
subpaths that are either single edges or exceptional subpaths. 
\end{enumerate}
\end{thm}
\qed

\begin{lem}\label{lem4}
Let $\mathcal{O} \in Out(F_n)$ be a UPG automorphism of infinite order 
and let  $f \co G \rightarrow G$ be its topological representative as in 
Theorem~\ref{thm2}. For every path $\gamma$ in $G$ for which 
$[[f(\gamma)]] \neq \gamma$ there exists $a \in \mathbb{R}$ such that
\[ \alpha([[f^k(\gamma)]]) \geq k+a \, .\]
\end{lem}

\begin{proof}
We prove our claim by induction on the (minimal) index, $m$, of the filtration 
element that contains a path $\gamma$. 

If $\gamma \subset G_1$ there is nothing to be proved since $G_1$ contains 
only one edge $E_1$ which is fixed by $f$. 

Suppose the claim is true for the subpaths contained in $G_{m-1}$ 
that satisfy out hypothesis, and let $\gamma$ be a path in $G_m$ for which  
$[[f(\gamma)]] \neq \gamma$. By Theorem~\ref{thm2} for every $m \geq 
M(\gamma)$, $[[f^m (\gamma)]]$ splits into subpaths that are either single 
edges or exceptional paths. 
Denote $[[f^{M(\gamma)} (\gamma)]]$ by $\tilde{\gamma}$, so that 
$\tilde{\gamma}=\gamma_1 \cdot \ldots \cdot \gamma_{p}$, where $\gamma_i$
is either a single edge or an exceptional path.

Assume there is an exceptional path $\gamma_t$ which is not fixed by $f$. 
Without loss of generality we may assume that $\gamma_t=E_i \upsilon^r 
\overline {E_j}$, where $f(E_i)=E_i \upsilon^l \  (l>0)$, $f(E_j)=E_j 
\upsilon^s \ (s>0)$ and $j \leq i$. Now we have that
\[ [[f^k(\gamma_t)]]=E_i \upsilon^{k(l-s)+r}\  \overline{E_j}\, , \]
and 
\[\alpha([[f^k(\gamma_t)]])\geq k(l-s)+r,\ \  \text{if} \ l-s>0\, ,\]
\[\alpha([[f^k(\gamma_t)]])\geq k(s-l)-r,\ \  \text{if}\ l-s<0\, .\]
Since $\gamma_t$ is not fixed, $l$ and $s$ cannot be equal. Therefore
\[\alpha([[f^k(\tilde{\gamma})]]) \geq k \pm r\, , \]
\[\alpha([[f^k(\gamma)]])=\alpha([[f^{k-M(\gamma)}(\tilde{\gamma})]]) \geq 
k -M(\gamma) \pm r\, . \]
If all exceptional paths in $\tilde{\gamma}$ are fixed, there exists an edge 
$\gamma_t=E_i$ which is not fixed by $f$. We know that $f(E_i)=E_i \cdot u_i$,
where $u_i$ is a closed path contained in $G_{m-1}$. \\ 
If $[[f(u_i)]] = u_i$, our claim is proven since $[[f^k(E_i)]]=E_i u_i^k$ 
and so 
\[\alpha([[f^k(\tilde{\gamma})]])\geq k \, ,\]
\[\alpha([[f^k(\gamma)]])\geq k-M(\gamma)\, . \]
If $[[f(u_i)]] \neq u_i$, there exists $a \in \mathbb{R}$ such that 
$\alpha([[f^k(u_i)]]) \geq k + a$. We now have 
\[\alpha([[f^k(\gamma)]])=\alpha([[f^{k-M(\gamma)}(\tilde{\gamma})]]) \geq 
\alpha([[f^{k-M(\gamma)}(u_i)]]) \geq k-M(\gamma)+a\, . \] 
\end{proof}

\begin{lem} \label{lem3}
Let $\mathcal{O}$ be a UPG automorphism of $F_n$ of infinite order. 
There exist a closed path $\sigma$ in $G$, and $b \in \mathbb{R}$  such that
\[\tilde{\alpha}(\mathcal{O}^k(\sigma))\geq k + b\, .\] 
\end{lem}

\begin{proof}
Let $f \co G \rightarrow G$ be as in Theorem~\ref{thm2}. Since $\mathcal{O} 
\neq id$ there is a closed path $\sigma$ which is not fixed by $f$.  We know 
that for every $m \geq M(\sigma)$, $[[f^m (\sigma)]]=\sigma_1 \cdot \ldots 
\cdot \sigma_{p}$ splits into subpaths that are either single edges or 
exceptional paths. 
Denote $[[f^{M(\sigma)} (\sigma)]]$ by $\tilde{\sigma}$, so that 
$\tilde{\sigma}=\sigma_1 \cdot \ldots \cdot \sigma_{p}$.

If there is an exceptional path  $\sigma_t$ in this splitting which is not 
fixed by $f$, we get 
\[\tilde{\alpha}(\mathcal{O}^k(\tilde{\sigma})) \geq k \pm r \]
as in Lemma~\ref{lem4}. 

If all exceptional paths in $\tilde{\sigma}$ are fixed, there exists an edge 
$\sigma_t=E_i$ such that $f(E_i)=E_i \cdot u_i$, where 
$u_i$ is a closed path contained in $G_{i-1}$. By Lemma~\ref{lem4} there 
exists $a \in \mathbb{R}$ such that
\[ \alpha(f^k(E_i)) \geq k+a\, .\]
Hence, in all the above cases, there is $b \in \mathbb{R}$ such that  
\[ \tilde{\alpha}(\mathcal{O}^k(\tilde{\sigma})) \geq k+b\, . \]
     
\end{proof}

\section{Proof of Theorem~\ref{thm1}}

We consider the cases of exponentially and polynomially growing outer 
automorphisms separately. 

{\bf Case 1.} Let $\mathcal{O}$ be an exponentially growing outer automorphism
of $F_n$. There exist $\lambda >1$ and a cyclically reduced word $w$ such 
that $\ell(\mathcal{O}^k([w])) \geq \lambda^{k} \ell([w])$, for all 
$k \geq 1$, where $\ell$ denotes the cyclic word length (see \cite{hb1}).  
Suppose that $\mathcal{O}^k$ can be written as $\tilde{g_1} \ldots 
\tilde{g_m}$, for some $\tilde{g_i} \in \tilde{Y}$. It is straightforward to 
show that for all $\tilde{g} \in \tilde{Y}$ and any cyclically reduced word 
$w$ we have 
\[ \ell (\tilde{g}([w])) \leq 2 \, \ell ([w]) \]
Using this inequality we obtain: 
\[ \lambda^{k} \ell ([w]) \leq \ell (\mathcal{O}^{k}([w])) \leq 2^{m} 
\ell ([w]) \]
Hence
\[ m \geq \frac{\log\lambda^{k}}{\log2} \]
which implies
\[\tau(\mathcal{O}) \geq \frac{\log\lambda}{\log2}>0\]\\
There is a constant $c_1 > 1$ such that $\lambda \geq c_1$ (\cite{hb1}). 
Therefore $\tau (\mathcal{O})$ is bounded away from zero. 

{\bf Case 2.} Let $\mathcal{O}$ be a $UPG$ automorphism. Again assume 
that $\mathcal{O}^k$ can be written as $\tilde{g_1}\ldots \tilde{g_m}$, 
for some $\tilde{g_i} \in \tilde{Y}$.  
By Lemma~\ref{lem3} there is a closed path $\sigma$ in $G$ such that    
\[\tilde{\alpha}(\mathcal{O}^k(\sigma))\geq k + b\] 
Let $u_j = \tilde{g_j} \ldots \tilde{g_m}$. Applying Lemma~\ref{lem1} we get 
\[ \tilde{\alpha}(u_i(\sigma)) \leq \tilde{\alpha}(u_{i+1}(\sigma)) + C \]
which yields
\begin{align*}
& k + b \leq \tilde{\alpha}(\mathcal{O}^k(\sigma)) \leq mC + 
\tilde{\alpha}(\sigma) \\
& \frac{k+b-\tilde{\alpha}(\sigma)}{C} \leq m \, .
\end{align*}
We have
\[ \tau(\mathcal{O}) \geq \lim_{k \rightarrow \infty} \frac{k+b-\tilde{\alpha}
(\sigma)}{k\, C} = \frac{1}{C}\, . \]
Finally, if $\mathcal{O}$ is any polynomially growing outer automorphism, 
then there exists $s \geq 1$ (bounded above by some $c_2$), such that 
$\mathcal{O}^s$ is a $UPG$ automorphism. Then
\[\tau(\mathcal{O})=\frac{1}{s}\, \tau(\mathcal{O}^s) \geq \frac{1}{C\, s} > 0
\, . \]
Since $s$ is bounded by $c_2$, we get $\tau (\mathcal{O}) \geq \frac{1}{C\, 
c_2} >0$.

This completes the proof. 
\qed

\bibliographystyle{halpha}
\bibliography{$GLOBAL/all}
\par

{\noindent \sc Department of Mathematics, University of Utah\\
155 S 1400 E, rm 233\\
Salt Lake City, UT 84112-0090, USA\\}
\par
{\noindent \it E-mail:} emina@math.utah.edu

\end{document}